\documentclass[12pt,oneside,american]{amsart}
\usepackage[T1]{fontenc}
\usepackage[utf8]{inputenc}
\usepackage{mathrsfs}
\usepackage{amstext}
\usepackage{amsthm}
\usepackage{amssymb}

\makeatletter
\numberwithin{equation}{section}
\numberwithin{figure}{section}
\theoremstyle{plain}
\newtheorem{thm}{\protect\theoremname}
\theoremstyle{plain}
\newtheorem{cor}{\protect\corollaryname}
\theoremstyle{plain}
\newtheorem{prop}{\protect\propositionname}
\theoremstyle{plain}
\newtheorem{lem}{\protect\lemmaname}

\usepackage{fullpage}

\makeatother

\usepackage{babel}
\providecommand{\corollaryname}{Corollary}
\providecommand{\lemmaname}{Lemma}
\providecommand{\propositionname}{Proposition}
\providecommand{\theoremname}{Theorem}

\begin{document}
\title{On the II$_{1}$ Factors of Fuchsian Groups}
\date{\today}
\address{Department of Mathematics, UCLA, Los Angeles, CA 90095}
\email{shlyakht@math.ucla.edu.}
\author{Dima Shlyakhtenko}
\begin{abstract}
We show that von Neumann algebras of fundamental groups of closed
orientable surfaces of genus $g\geq2$ are free group factors on $2g-1$
generators. The key technical ingredient involves a proof that the
element $w=ABA^{-1}B^{-1}$ of the free group $\mathbb{F}_{2}=\langle A,B\rangle$
is freely complemented in the group factor: $L(\mathbb{F}_{2})=W^{*}(w)*W^{*}(v)$
for some Haar unitary $v\in L(\mathbb{F}_{2})$ that is freely independent
from $w$. Combined with previous results, we conclude that for an
arbitrary finitely generated torsion-free non-elementary discrete
subgroup $\Gamma\subset PSL_{2}(\mathbb{R})$, $L(\Gamma)$ is a free
group factor, settling a conjecture of de la Harpe and Voiculescu.
This result was obtained using OpenAI's ChatGPT Pro 6.0.
\end{abstract}

\maketitle

\section{Introduction.}

Recall (see e.g. \cite{delaHarpe-Voiculescu:fuchsian}) that a finitely
generated non-elementary (i.e., not finite and not a finite extension
of $\mathbb{Z}$) discrete subgroup of $PSL_{2}(\mathbb{R})$ has
the presentation

\begin{align*}
\Gamma(g:m_{1},\dots,m_{k}:t,s) & =\Bigg\langle a_{1},b_{1},\ldots,a_{g},b_{g},x_{1},\ldots,x_{k},p_{1},\dots,p_{s},h_{1},\dots,h_{t}\ \ \Big|\\
 & \forall j,x^{m_{j}}=1,\quad\left(\prod_{i=1}^{g}a_{i}b_{i}a_{i}^{-1}b_{i}^{-1}\right)\cdot x_{1}\cdots x_{k}\cdot p_{1}\cdots p_{s}\cdot h_{1}\cdots h_{t}=1\Bigg\rangle
\end{align*}
where $g,k,t,s\geq0$, $m_{1},\dots,m_{k}\geq2$ are satisfy the condition
\[
r:=2g-1+s+t+\sum_{j=1}^{r}\left(1-\frac{1}{m_{j}}\right)>1.
\]

In \cite[Problem on p. 157]{delaHarpe-Voiculescu:fuchsian} de la
Harpe and Voiculescu asked whether it's always true that $L(\Gamma(g:m_{1},\dots,m_{k},t,s))\cong L(\mathbb{F}_{r})$;
moreover, they showed that this is true if $t+s>0$. (If $r$ is non-integer,
$L(\mathbb{F}_{r})$ stands for an interpolated free group factor
\cite{dykema:interpolated,radulescu:interpolated}.)

If the group is torsion-free (i. e. $k=0$) and $t+s>0$, the the
answer is immediate (as observed in \cite{delaHarpe-Voiculescu:fuchsian}),
since the group itself is a free group. The remaining torsion free
case corresponds to $t=s=k=0$ and $g\geq2$. In that case, the group
is naturally an amalgamated free product of two free groups, with
amalgamation over an infinite cyclic subgroup:
\[
\Gamma(g::0,0)=\langle a_{1},b_{1}\rangle*_{a_{1}b_{1}a_{1}^{-1}b_{1}^{-1}=(\prod_{j=2}^{g}a_{j}b_{j}a_{j}^{-1}b_{j}^{-1})^{-1}}\langle a_{2},b_{2},\dots,a_{g},b_{g}\rangle
\]
and its isomorphism class remained open until now (see the discussion
of the Problem on p. 157 of \cite{delaHarpe-Voiculescu:fuchsian}). 

The main technical result of this paper, Theorem \ref{thm:freeComplementation},
shows that there exists a Haar unitary $v$, for which $L(\mathbb{F}_{2})=W^{*}(a_{1},b_{1})=W^{*}(a_{1}b_{1}a_{1}^{-1}b_{1}^{-1})*W(v)$.
In other words, the abelian von Neumann subalgebra of $L(\mathbb{F}_{2})$
generated by the commutator $a_{1}b_{1}a_{1}^{-1}b_{1}^{-1}$ is freely
complemented (see \cite[Remark 1.4(2)]{Boutonnet-Popa:freeComplementation}
for a discussion of this question; see also \cite[\textbackslash{}S3.2]{popa:compressible}
and \cite{popa:maxInjective}). This implies (see Corollary \ref{cor:amalgamatedIsFree})
that the amalgamated free product above is indeed a free group factor:
\[
W^{*}(v)*W^{*}(a_{2},b_{2},\dots,a_{g},b_{g})\cong L(\mathbb{F}_{2g-1}).
\]

The following theorem summarizes all known results for discrete finitely
generated non-elementary subgroups of $PSL_{2}(\mathbb{R})$; cases
2 and 3 are new and are proved in this paper:
\begin{thm}
\label{thm:subgroupsPSL2}Suppose that any of the following conditions
hold:
\begin{enumerate}
\item $t+s>0$
\item $g\geq2$
\item $g\geq1$ and $k\geq2$
\end{enumerate}
Then $L(\Gamma(g;m_{1},\dots,m_{k}:t,s)\cong L(\mathbb{F}_{r})$ where
$r=2g-1+s+t+\sum_{j=1}^{q}(1-\frac{1}{m_{j}})$. 
\end{thm}
\begin{proof}
The first case corresponds to cases (2) and (4) described on p. 156
of \cite{delaHarpe-Voiculescu:fuchsian}. As noted there, the group
$\Gamma(g:m_{1},\dots,m_{k}:t,s)$ is a free product of cyclic groups,
and $L(\Gamma(g;m_{1},\dots,m_{k}:t,s)\cong L(\mathbb{F}_{r})$ by
the results of Dykema \cite{dykema:freeProductHyperfinite}.

In the second and third cases, consider the group $K$ freely generated
by $a_{2},b_{2},\dots,a_{g},b_{g}$ as well as $x_{1},\dots,x_{k}$
subject only to $x_{j}^{m_{j}}=e$ (if $g=1$, we only have $x_{1},\dots,x_{k}$
as generators). Thus
\[
K=\mathbb{Z}/m_{1}\mathbb{Z}*\cdots*\mathbb{Z}/m_{k}\mathbb{Z}*(\mathbb{Z})^{*(2g-2)}.
\]

Put 
\[
z=\prod_{i=2}^{g}a_{i}b_{i}a_{i}^{-1}b_{i}^{-1}\cdot x_{1}\cdots x_{k}
\]
Then the assumptions that $g\geq2$ or $k\geq2$ guarantee that $z$
has infinite order. 

Let $\Lambda=\mathbb{F}_{2}$ generated by two elements $a_{1}$ and
$b_{1}$. Then
\[
\Gamma(g:m_{1},\dots,m_{k}:0,0)\cong\Lambda*_{a_{1}b_{1}a_{1}^{-1}b_{1}^{-1}=z^{-1}}K.
\]
As before, let $\mathscr{A}=W^{*}(aba^{-1}b^{-1})\subset L(\Lambda)$.
Applying Corollary \ref{cor:amalgamatedIsFree}
\begin{align*}
L(\Gamma(g & :m_{1},\dots,m_{k}:0,0))\cong L(\mathbb{F}_{2})*_{\mathscr{A}}L(K)\\
 & \cong L(\mathbb{Z})*L(K)\cong L(\mathbb{Z}*K)\\
 & \cong L(\mathbb{Z}/m_{1}\mathbb{Z}*\cdots*\mathbb{Z}/m_{k}\mathbb{Z}*(\mathbb{Z})^{*(2g-1)})\\
 & \cong L(\Gamma(0:m_{1},\dots,m_{k}:2g)\cong L(\mathbb{F}_{2g-1+\sum(1-m_{j}^{-1})})
\end{align*}
where the last isomorphism uses the first case of the Theorem.
\end{proof}
\begin{cor}
\label{cor:torsionFreeAreFree}If $\Lambda$ is a torsion-free finitely
degenerated non-elementary discrete subgroup of $PSL_{2}(\mathbb{R})$,
then $L(\Gamma)\cong L(\mathbb{F}_{r})$ (where $r$ is as above).
In particular, if $\Lambda=\pi_{1}(\Sigma_{g})$ is the fundamental
group of a surface of genus $g\geq2$, then $L(\Lambda)\cong L(\mathbb{F}_{2g-1}).$
\end{cor}
\begin{proof}
The torsion-free case corresponds to $k=0$. Since $r=1+2g-2+s+t>1$
this means that either $s+t>0$, or $g\geq2$, corresponding to parts
1 or 2 of Theorem \ref{thm:subgroupsPSL2}.
\end{proof}
While we are unable to settle the de la Harpe-Voiculescu conjecture
in all cases, we are able to prove a weaker statement of containment
in (or of) an interpolated free group factor up to finite index:
\begin{thm}
Let $\Gamma$ be a finitely generated non-elementary discrete subgroup
of $PSL_{2}(\mathbb{R})$. Then $L(\Gamma)$ contains a free group
factor as a finite-index subfactor. Moreover, $L(\Gamma)$ is contained
in an interpolated free group factor as a finite-index subfactor.
\end{thm}
\begin{proof}
The remaining families of discrete subgroups of $PSL_{2}(\mathbb{R})$
not covered by Theorem \ref{thm:subgroupsPSL2} are:

The cases $\Gamma(0:m_{1},\dots,m_{k}:0,0)$ with $m_{1},\dots,m_{k}$
satisfying $\sum(1-m_{j}^{-1})>2$ (thus necessarily $k\geq3$). The
group presentation is $\left\langle x_{1},\ldots,x_{k}\ |\ x_{j}^{m_{j}}=1,\ x_{1}\cdots x_{k}=1\right\rangle $.
The prediction is that the group von Neumann algebra should be $L(\mathbb{F}_{r})$
with $r=(\sum(1-m_{j}^{-1}))-1$.

The cases $\Gamma(1:m:0,0)$. The presentation is $\left\langle a,b\ |\ (aba^{-1}b^{-1})^{m}=1\right\rangle $,
for some $m\geq2$. The prediction is that the von Neumann algebra
should be $L(\mathbb{F}_{r}$) with $r=2-1/m$.

In either of these cases, $\Gamma$ is a co-compact subgroup of $PSL_{2}(\mathbb{R})$.
By the solution of the Fenchel's Conjecture (see \cite{nielsen,Bundgaard-Nielsen,fox,Chau,Mennicke,Mennicke:correction}),
there exists a finite-index torsion-free subgroup $\Lambda\subset\Gamma$.
Since $\Gamma\subset PSL_{2}(\mathbb{R})$ is co-compact, $\Lambda$
is also co-compact and thus by Corollary \ref{cor:torsionFreeAreFree},
$L(\Lambda)\cong L(\mathbb{F}_{k})$ for some integer $k$. Thus $L(\mathbb{F}_{k})\cong L(\Lambda)\subset L(\Gamma)$
has finite index. 

Using the Jones basic construction \cite{jones:index} we also obtain
the inclusion $L(\Gamma)\subset\langle L(\Gamma),e_{L(\Lambda)}\rangle$.
By \cite{pimsner-popa:index}, $\langle L(\Gamma),e_{L(\Lambda)}\rangle$
is the amplification $L(\Lambda)^{[\Gamma:\Lambda]}$. Using the amplification
formula for free group factors \cite{dykema:interpolated,dykema:freeProductHyperfinite,radulescu:interpolated}
we conclude that $M_{1}$ is also an intorpolated free group factor.
\end{proof}
It is worth mentioning that our result can be viewed as a group von
Neumann algebra counterpart to the work of Conley, Gaboriau, Marks
and Tucker-Drob (\cite{conleyGaboriauMarksTuckerDrob:treeabilitySurface},
see also \cite{gaboriau:stronglyTreeable}) where a similar ``free
complementation`` result was proved in the context of equivalence
relations generated by group actions with planar Cayley graphs (in
particular, fundamental groups of surfaces).

\subsection{How this result was obtained: use of AI}

Our original approach to the free complementation question (from some
years ago) was to construct a free complement to $w$ by finding a
family of ``changes of generators'' that would transform $A,B$
into the pair $(w,v)$ and thus would yield a free complement $v$
to $w$, so that $W^{*}(u,v)$ still generates all of $W^{*}(A,B)$.
This approach is natural from free transportation theory. It is also
natural to try to create such a path of automorphisms $\alpha_{t}$
by solving a differential equation, in effect, integrating a time-dependent
divergence-free vector field. To keep the changes of variables invertible,
one needs to also set up a differential equation for the inverse transformation,
and to show that these differential equations have solutions that
behave well. Since not every unitary in $L(\mathbb{F}_{2})$ is freely
complemented (or even generates a MASA) one needs to extract specific
analytical features from the particular case of $w=ABA^{-1}B^{-1}$,
and our initial starting point was certain combinatorial considerations
connected to the proof of the ergodic theory version of the problem
\cite{conleyGaboriauMarksTuckerDrob:treeabilitySurface}. ChatGPT
Pro 5.6 was then used to experiment with various approaches: showing
existence of solutions, invertibility of changes of variables, etc.,
all aiming to obtain and propagate estimates involving the particular
geometry of the subgroup generated by $w$. The manuscript was approaching
70 pages, but the we were still unable to find the right estimates
that would close the proof. We then asked ChatGPT Pro 6.0 to take
a look at the manuscript and suggest improvements in the approach.
After some thought, it returned a 7-page proof that became the basis
of this manuscript. In its own words:
\begin{quotation}
The argument arose from an attempt to obtain free complementation
by integrating constrained time-dependent derivations, motivated by
measured free-factor constructions of \cite{conleyGaboriauMarksTuckerDrob:treeabilitySurface}.
Although the proof presented here does not use that transport construction,
the elementary commutator-preserving automorphisms and the need to
preserve generation in the limit emerged from its analysis. The resulting
argument instead uses finite changes of generators, the associated
tree geometry, and quantitative preservation of polynomial approximations.

This work was developed with substantial assistance from ChatGPT (OpenAI),
used interactively to explore proof strategies, propose and critique
intermediate arguments, carry out exploratory calculations, and draft
and revise the manuscript. This assistance included mathematical content,
in particular the development of the finite-automorphism approach
and its Fourier estimates, and was not limited to editing. The author
supplied the research question and initial geometric motivation, and
directed the investigation.

The development proceeded through an extended sequence of proposals,
objections, corrections, and changes of approach; the final argument
should not be understood as an unmodified output from a single prompt.
\end{quotation}
The argument presented by ChatGPT was then carefully checked and rewritten
by the author. The author is fully responsible for the content of
this manuscript, including the mathematical claims, proofs, and references.

\subsection{Outline of the construction. }

The free complement $v$ is obtained as a limit of elements $v_{n}$,
which are themselves constructed as $v_{n}=(\beta_{1}\circ\cdots\circ\beta_{n})(A)$,
where $\beta_{1},\dots,\beta_{n}$ are certain automorphisms of $L(\mathbb{F}_{2})$
that satisfy $\beta_{j}(w)=w$. Since $A$ and $w$ are freely independent,
it follows that $v_{n}$ and $w$ are free at all times. The automorphisms
$\beta_{j}$ are chosen in such a way that $v_{n}$ form a Cauchy
sequence, so that $v_{n}\to v$. What is key, however, is to ensure
that $W^{*}(v,w)=L(\mathbb{F}_{2})$. This is guaranteed by choosing
automorphisms in such a way that both $A$ and $B$ are well-approximated
by $v_{n}$ and $w$. Writing $N=W^{*}(A,w)$, the approximation is
achieved by choosing $\Theta_{n}$ so that the conjugated conditional
expectation $E_{\Theta_{n}(N)}=\Theta_{n}\circ E_{N}\circ\Theta_{n}^{-1}$
approximates the identity on prescribed finite sets. Proposition \ref{prop:ExistsAutomorphism}
formalizes the passage from such local approximation to free complementation
through a Baire category argument. In proposition \ref{prop:TandTheta}
we show that instead of controlling $\Theta_{n}\circ E_{W^{*}(A,w)}\circ\Theta_{n}^{-1}$,
one can instead use a family of positive $L^{2}$ contractions $\mathbf{T}_{k}$
interpolating between $E_{W^{*}(A,w)}$ and identity, and prove a
certain inequality involving the quadratic form $\mathbf{Q}_{k}$
associated to $\mathbf{T}_{k}$. 

The remainder of the paper is dedicated to producing a family of automorphisms
that satisfy the hypothesis of Proposition \ref{prop:TandTheta} and
the choice of the family $\mathbf{T}_{k}$. The automorphisms $\beta_{j}$
are obtained by repeatedly composing automorphism of the form
\[
R_{-m}\circ S_{f_{m,t}}\circ R_{m}
\]
where $R_{m}(A)=A$, $R_{m}(B)=BA^{m}$, and $S_{f}(B)=B$, $S_{f_{m,t}}(A)=Af(B)$
for certain functions $f_{m,t}:\mathbb{T}\to\mathbb{T}$ . The contractions
$\mathbf{T}_{k}$ are diagonal in the group basis, and multiply an
element $h$ by $q_{k}^{\ell(h)}$ where $\ell(h)$ is the length
of $h$ relative to the subgroup generated by $w$ and $A$, and $q_{k}\uparrow1$,
$q_{0}=0$ are carefully chosen, along with a family of functions
$f_{m,t}$. The key estimate is given by Lemma \ref{lem:F_tlowerEqThetaOfB}.

\subsection*{Acknowledgements.}

This research was sponsored in part by the Army Research Office and
was accomplished under Grant Number W911NF-25-1-0075. The views and
conclusions contained in this document are those of the authors and
should not be interpreted as representing the official policies, either
expressed or implied, of the Army Research Office or the U.S. Government.
The U.S. Government is authorized to reproduce and distribute reprints
for Government purposes notwithstanding any copyright notation herein.
Research was also supported in part by NSF grant DMS-2348633.

\section{Free complementation.}

Let $\mathbb{F}_{2}=\langle A,B\rangle$ be the free group on two
generators, and let $M=L(\mathbb{F}_{2})$. Regard $\mathbb{F}_{2}\subset M$.
We endow $M$ with its unique normalized group trace $\tau$. 

Suppose that $w\in M$ is a Haar unitary so that $w$ is free from
$A$, and put $N=W^{*}(w,A)\subset M$. We give the following sufficient
condition for the existence of an abelian free complement to $w$,
i.e., the existence of a Haar unitary $v$ free from $w$ and so that
$W^{*}(w,v)=M$:
\begin{prop}
\label{prop:ExistsAutomorphism}Suppose that for any $x_{1},\dots,x_{n}\in M$
and any $\epsilon>0$, $\eta>0$, there exists an automorphism $\theta:M\to M$
so that
\begin{itemize}
\item $\theta(w)=w$
\item $\Vert\theta(A)-A\Vert_{2}<\eta$
\item $\sum_{j=1}^{n}\Vert x_{j}-E_{\theta(N)}(x_{j})\Vert_{2}^{2}<\epsilon$. 
\end{itemize}
Then for any $\delta>0$ there exists a Haar unitary $v\in M$ so
that:
\begin{itemize}
\item $v$ is free from $w$, 
\item $M=W^{*}(w,v)=W^{*}(w)*W^{*}(v)$ 
\item $\Vert v-A\Vert_{2}<\delta$.
\end{itemize}
In particular, $W^{*}(w)$ is freely complemented in $M$ by the abelian
algebra $W^{*}(v)$.
\end{prop}
\begin{proof}
Let $X=\overline{\{\theta(A):\theta\in\operatorname{Aut}(M),\ \theta(w)=w\}}^{L^{2}(M)}$.
Since the unitary group of $M$ is closed in the $L^{2}$ norm and
$X$ is a closed subset, $X$ is complete. Since $A$ is a Haar unitary
free from $w$, any automorphism fixing $w$ takes $A$ to a Haar
unitary free from $w$. Since freeness is preserved under limits,
it follows that every element of $X$ is a Haar unitary free from
$w$. 

Now choose a total sequence $x_{k}\in M$ for $\Vert\cdot\Vert_{2}$.
For each $n>0$ let
\[
O_{n}=\{u\in X:\exists\textrm{ \ensuremath{*}-polynomials \ensuremath{P_{1},\dots,P_{n}} with }\max_{1\leq j\leq n}\Vert x_{j}-P_{j}(w,u)\Vert_{2}<1/n\}.
\]
Since each $P_{j}$ is Lipschitz in its unitary arguments for the
$L^{2}$ norm, the set $O_{n}$ is open. Furthermore, $O_{n}$ is
dense. Indeed, suppose that $\theta_{0}\in\operatorname{Aut}(M)$,
$\theta_{0}(w)=w$, and $\rho>0$ are given. Then we can find an automorphism
$\theta$ fixing $w$ and such that $\Vert\theta(A)-A\Vert_{2}<\rho$
and
\[
\sum_{j=1}^{n}\Vert\theta_{0}^{-1}(x_{j})-E_{\theta(N)}(\theta_{0}^{-1}(x_{j}))\Vert_{2}^{2}<1/(4n^{2}).
\]
Equivalently,
\[
\sum_{j=1}^{n}\Vert x_{j}-E_{\theta_{0}(\theta(N))}(x_{j})\Vert_{2}^{2}<1/(4n^{2}).
\]
Thus $x_{j}$ is approximated by elements $E_{\theta_{0}(\theta(N))}(x_{j})$
of $\theta_{0}(\theta(N))=W^{*}(\theta_{0}(\theta(A)),w)$ to within
$1/(2n)$ in $L^{2}$. Therefore, we can find $*$-polynomials $P_{1},\dots,P_{n}$
with the property that $P_{j}(w,\theta_{0}(\theta(A)))$ are within
$1/2n$ of $E_{\theta_{0}(\theta(N))}(x_{j})$ and so within $1/n$
of $x_{j}$. Thus $(\theta_{0}\circ\theta)(A)\in O_{n}$. On the other
hand, $\Vert\theta_{0}(\theta(A))-\theta_{0}(A)\Vert_{2}<\rho$ so
that $\theta_{0}(A)$ is within $\rho$ of $O_{n}$. Since the the
set $\{\theta(A):\theta\in\operatorname{Aut}(M),\ \theta(w)=w\}$
is dense in $X$, it follows that $O_{n}$ is dense as well. 

By Baire Category Theorem, the intersection of the sets $O_{n}$ is
dense. Choose a unitary $v$ in this intersection so that $\Vert v-A\Vert_{2}<\delta$.
Then $v$ is free from $w$, and, since for each $j$ , $v\in O_{n}$
for arbitary large $n>j$, it follows that $x_{j}\in L^{2}(W^{*}(w,v))$.
Since $\{x_{j}\}_{j}$ is total, it follows that $W^{*}(w,v)=M$.
\end{proof}
\begin{prop}
\label{prop:TandTheta}Suppose that there exist positive contractions
$\mathbf{T}_{k}:L^{2}(M)\to L^{2}(M)$, $k\geq0$, so that $\mathbf{T}_{0}=E_{N}$
and $\mathbf{T}_{k}\to\operatorname{Id}$ strongly. Put $\mathbf{Q}_{j}(x)=\langle\mathbf{T}_{j}x,x\rangle$. 

Suppose moreover that for any $k,n\geq0$ and any $\epsilon>0,\eta>0$,
and any $x_{1},\dots,x_{n}\in M$, there exists an automorphism $\theta:M\to M$
so that:
\begin{itemize}
\item $\theta(w)=w$
\item $\Vert\theta(A)-A\Vert_{2}<\eta$
\item $\sum_{i=1}^{n}\mathbf{Q}_{k}(\theta^{-1}(x_{i}))\geq\sum_{i=1}^{n}\mathbf{Q}_{k+1}(x_{i})-\epsilon$.
\end{itemize}
Then the hypothesis of Proposition \ref{prop:ExistsAutomorphism}
is satisfied, and thus its conclusion holds.
\end{prop}
\begin{proof}
Let $x_{1},\dots,x_{n}$, $\eta,\epsilon$ be given as in the hypothesis
of Proposition \ref{prop:ExistsAutomorphism}. Choose $J\geq1$ large
enough so that 
\[
\sum_{i}\mathbf{Q}_{J}(x_{i})>\sum_{i}\Vert x_{i}\Vert_{2}^{2}-\epsilon/2.
\]
Let $x_{i}^{(J)}=x_{i}$, and construct $x_{i}^{(J-1)},\dots,x_{i}^{(0)}$
recursively as follows. Having chosen $x_{i}^{(j)}$ and automorphisms
$\theta_{J-1},\dots,\theta_{j}$, choose $\theta_{j-1}$ so that $\theta_{j-1}(w)=w$,
$\Vert\theta_{j-1}(A)-A\Vert_{2}<\eta/2J$,
\[
\sum_{i}\mathbf{Q}_{j-1}(\theta_{j-1}^{-1}(x_{i}^{(j)}))\geq\sum_{i}\mathbf{Q}_{j}(x_{i}^{(j)})-\epsilon/2J,
\]
and set $x_{i}^{(j-1)}=\theta_{j-1}^{-1}(x_{i}^{(j)})$. Then
\[
\sum_{i}\mathbf{Q}_{j-1}(x_{i}^{(j-1)})\geq\sum_{i}\mathbf{Q}_{j}(x_{i}^{(j)})-\epsilon/2J.
\]

Let now $\Theta=\theta_{J-1}\circ\cdots\circ\theta_{0}$. Since $\mathbf{T}_{0}=E_{N}$
is an orthogonal projection, and noting that $\Theta^{-1}(x_{i})=x_{i}^{(0)}$,
we get that $\Vert\Theta(A)-A\Vert_{2}<\eta$ and
\begin{align*}
\sum_{i}\Vert x_{i}-E_{\Theta(N)}(x_{i})\Vert_{2}^{2} & =\sum_{i}\Vert\Theta^{-1}(x_{i})-E_{N}(\Theta^{-1}(x_{i}))\Vert_{2}^{2}\\
 & =\sum_{i}\left[\Vert x_{i}\Vert_{2}^{2}-\langle E_{N}(x_{i}^{(0)}),x_{i}^{(0)}\rangle\right]\\
 & =\sum_{i}\left[\Vert x_{i}\Vert_{2}^{2}-\mathbf{Q}_{0}(x_{i}^{(0)})\right]\\
 & \leq\sum_{i}\left[\Vert x_{i}\Vert_{2}^{2}-\mathbf{Q}_{1}(x_{i}^{(1)})\right]+\epsilon/2J\\
 & \leq\cdots\leq\sum_{i}\left[\Vert x_{i}\Vert_{2}^{2}-\mathbf{Q}_{J}(x_{i}^{(J)})\right]+\epsilon/2\\
 & <\sum_{i}\left[\Vert x_{i}\Vert_{2}^{2}-\Vert x_{i}\Vert_{2}^{2}\right]+\epsilon/2+\epsilon/2=\epsilon,
\end{align*}
Similarly, $\Vert\Theta(A)-A\Vert_{2}\leq\sum_{j=0}^{J-1}\Vert\theta_{j}(A)-A\Vert_{2}\leq\eta/2<\eta$. 
\end{proof}

\section{The case of the commutator MASA.}

Let as before $M=L(\mathbb{F}_{2})$ with group generators $A$ and
$B$, and let $w=ABA^{-1}B^{-1}$ be their group commutator. Note
that $w$ and $A$ are free inside $\mathbb{F}_{2}$. Let $H$ denote
the subgroup generated by $A$ and $w$ . Set $N=W^{*}(A,w)=L(H)$. 

\subsection{A class of automorphisms.}

We first describe a class of automorphisms of $M$ that will be used
to deduce the hypothesis of Proposition \ref{prop:TandTheta}.
\begin{lem}
Let $f:\mathbb{T}\to\mathbb{T}$ be a measurable function. Then the
map $A\mapsto Af(B)$, $B\mapsto B$ extends to a normal trace-preserving
automorphism $S_{f}:M\to M$. The inverse of $S_{f}$ is $S_{\bar{f}}$
, and $S_{f}(w)=w$. 
\end{lem}
\begin{proof}
Since $A$ and $B$ are two free Haar unitaries, it's immediate to
verify that also $Af(B)$ and $B$ are free Haar unitaries. Thus the
map $A\mapsto Af(B)$, $B\to B$ extends to a normal injective $*$-homomorphism
$S_{f}$ from $M$ to itself. Since $B\in W^{*}(Af(B),B)$, also $A=(Af(B))\bar{f}(B)\in W^{*}(Af(B),B)$
and thus $W^{*}(Af(B),B)=M$. Hence $S_{f}$ is onto and thus an automorphism.
Clearly, $S_{\bar{f}}S_{f}(B)=B$ and $S_{\bar{f}}S_{f}(A)=S_{\bar{f}}(Af(B))=A\bar{f}(B)f(B)=A$,
so $S_{\bar{f}}$ is the inverse to $S_{f}$. Finally, $S_{f}(w)=S_{f}(ABA^{-1}B^{-1})=Af(B)Bf(B)^{-1}A^{-1}B^{-1}=w.$
\end{proof}
\begin{lem}
For each $m\in\mathbb{Z}$, the map $A\mapsto A$, $B\mapsto BA^{m}$
extends to a normal trace-preserving automorphism $R_{m}:M\to M$.
The inverse of $R_{m}$ is $R_{-m}$, and $R_{m}(w)=w$. 
\end{lem}
\begin{proof}
The automorphism $R_{m}$ is conjugate to $S_{z^{m}}$ via the map
that exchanges $A$ and $B$. Moreover, $R_{m}(w)=ABA^{m}A^{-1}A^{-m}B^{-1}=w$. 
\end{proof}
It will be useful to give a more detailed description of the automorphism
$S_{f}$. To this end, suppose that $g\in\mathbb{F}_{2}$ is fixed.
Suppose that $g$, as a reduced word in generators $A$ and $B$,
contains a total of $p\geq1$ occurrences of $A^{\pm1}$ (if $p=0$,
$S_{f}(B^{r})=B^{r})$. Then there exists a unique $(2p+1)$-tuple
$(\epsilon_{1},\dots,\epsilon_{p},n_{0},\dots,n_{p})$, where $\epsilon_{j}\in\{-1,1\}$,
$n_{j}\in\mathbb{Z}$ and $n_{j}\neq0$ whenever $\epsilon_{j}=-\epsilon_{j+1}$
so that
\[
g=B^{n_{0}}A^{\epsilon_{1}}B^{n_{1}}A^{\epsilon_{2}}\cdots A^{\epsilon_{p}}B^{n_{p}}\in\mathbb{F}_{2}.
\]
For example, if $g=A^{2}$, the associated tuple is $(1,1,0,0,0)$;
and if $g=ABA^{-1}$, the associated tuple is $(1,-1,0,1,0).$

Let $d_{0}=-1$ if $\epsilon_{1}=-1$ and $d_{0}=0$ otherwise, $d_{p}=1$
if $\epsilon_{p}=1$ and $d_{p}=0$ otherwise; and let
\[
d_{j}=\frac{\epsilon_{j}+\epsilon_{j+1}}{2},\qquad1\leq j<p.
\]

Let now $f:\mathbb{T}\to\mathbb{T}$ be measurable, and put $c_{j}(k)=\widehat{f^{d_{j}}}(k)=\int f(z)^{d_{j}}z^{-k}d\mu(z)$,
where $\mu$ is the normalized Haar measure. Let $K$ be the set of
$p+1$-tuples $(k_{0},\dots,k_{p})$ with $k_{j}\in\mathbb{Z}$ and
such that $k_{j}=0$ whenever $d_{j}=0$. Set
\[
W(k_{0},\dots,k_{p})=B^{n_{0}+k_{0}}A^{\epsilon_{1}}B^{n_{1}+k_{1}}A^{\epsilon_{2}}\cdots A^{\epsilon_{p}}B^{n_{p}+k_{p}}.
\]
The map $W$ is injective. Indeed, the only cancellation possible
would be if for some $j$, $\epsilon_{j}\neq\epsilon_{j+1}$, so that
$d_{j}=0$. But in that case $k_{j}=0$ and so $n_{j}+k_{j}=n_{j}\neq0$. 
\begin{lem}
\label{lem:descriptionOfSfWords}Let $g\in\mathbb{F}_{2}$ be a reduced
word containing $p\geq1$ occurrences of $A^{\pm1}$, $f:\mathbb{\mathbb{T}\to\mathbb{T}}$
measurable, and let $c_{j}(k)$, $W$, and $K$ be as above. Then:
\begin{itemize}
\item If $h$ is not in the range of the map $W$, then $\langle h,S_{f}(g)\rangle=0.$
\item If $h=W(k_{0},\dots,k_{p})$, then $\langle h,S_{f}(g)\rangle=\prod_{j=0}^{p}c_{j}(k_{j}).$
\end{itemize}
\end{lem}
\begin{proof}
By definition of $S_{f}$ and $d_{j}$,
\[
S_{f}(g)=B^{n_{0}}f^{d_{0}}(B)A^{\epsilon_{1}}B^{n_{1}}f^{d_{1}}(B)\cdots A^{\epsilon_{p}}B^{n_{p}}f^{d_{p}}(B).
\]
Substituting $f^{d_{j}}(B)=\sum c_{j}(k)B^{k}$ and expanding writes
$S_{f}(g)$ as the orthogonal $L^{2}$ sum of words $W(k_{0},\dots,k_{p})$
each occurring with the coefficient $\prod_{j=0}^{p}c_{j}(k_{j})$. 
\end{proof}
\begin{lem}
\label{lem:OffDiagonalEstimate}Let $f:\mathbb{T}\to\mathbb{T}$ as
before. Let $D$ be a contraction which is diagonal in the group basis
for $L^{2}(M)$, i.e. $Dg=\lambda_{g}g$, $|\lambda_{g}|\leq1$. Then
for all $g,h\in\mathbb{F}_{2}$, $g\neq h$,
\[
\left|\langle S_{f}(g),DS_{f}(h)\rangle\right|\le2\sqrt{\eta_{f}}+\eta_{f},
\]
where $\eta_{f}=1-|\tau(f(B))|^{2}.$
\end{lem}
\begin{proof}
Suppose first that $h$ has no occurrences of $A$, so $h=B^{r}$.
Then $DS_{f}(h)=\lambda_{h}h$. If $g$ has at least one occurrence
of $A^{\pm1}$, Lemma \ref{lem:descriptionOfSfWords} implies that
$S_{f}(g)\perp h$ and so the inner product is zero. For the same
reason, the inner product is also zero if $h$ has at least one occurrence
of $A^{\pm1}$ and $g$ does not. If neither has any occurrence of
$A^{\pm1}$, then $g\ne h$ implies that $S_{f}(g)=g\perp\lambda_{h}h=DS_{f}(h)$.
Thus the only interesting case is when both $g$ and $h$ have at
least one occurrence of $A^{\pm1}$. 

If this is the case, then Lemma \ref{lem:descriptionOfSfWords} states
that for certain integers $p$ and $p'$ , $n_{0},\dots,n_{p},n_{0}',\dots,n_{p'}'$
signs $\epsilon_{1},\dots,\epsilon_{p}$, $\epsilon'_{1},\dots\epsilon'_{p'}$,
and sets $K,K'$ we have:
\begin{align*}
S_{f}(g) & =\sum_{(k_{0},\dots k_{p})\in K}\left(\prod_{j=0}^{p}c_{j}(k_{j})\right)B^{n_{0}+k_{0}}A^{\epsilon_{1}}B^{n_{1}+k_{1}}A^{\epsilon_{2}}\cdots A^{\epsilon_{p}}B^{n_{p}+k_{p}}\\
S_{f}(h) & =\sum_{(k'_{0},\dots k'_{p'})\in K'}\left(\prod_{j=0}^{p'}c_{j}'(k'_{j})\right)B^{n_{0}'+k_{0}'}A^{\epsilon'_{1}}B^{n'_{1}+k'_{1}}A^{\epsilon'_{2}}\cdots A^{\epsilon'_{p'}}B^{n'_{p'}+k'_{p'}}.
\end{align*}
Thus $S_{f}(g)$ and $S_{f}(h)$ are linear combinations of disjoint
sets of words (and thus $\langle S_{f}(g),DS_{f}(h)\rangle=0$ since
$D$ is diagonal) unless 
\[
p=p'\quad\text{and}\quad\epsilon_{j}=\epsilon'_{j}\qquad(1\le j\le p).
\]
We may therefore assume these equalities hold. In particular, 
\[
d_{j}=d'_{j},\qquad c_{j}=c'_{j}\qquad(0\le j\le p).
\]
If $d_{j}=0$, then $k_{j}=k'_{j}=0$ by definition of $K$ and $K'$.
Thus a word can occur in both sums only if
\[
n_{j}+k_{j}=n'_{j}+k'_{j},
\]
which implies $n_{j}=n'_{j}$. Consequently, the inner product is
again zero if $n_{j}\ne n'_{j}$ whenever $d_{j}=0$. Thus we may
further assume that $n_{j}=n'_{j}$ whenever $d_{j}=0$. Put 
\[
r_{j}=n'_{j}-n_{j}.
\]
Since $g\ne h$ and their sign sequences agree, at least one $r_{j}$
is nonzero. Let us assume that $r_{m}\neq0$ (and thus $d_{m}\neq0$). 

Expanding the sums in the expressions for $S_{f}(g)$ and $S_{f}(h)$,
we see that the only terms that contribute to the inner product are
ones where $n_{j}+k_{j}=n'_{j}+k'_{j}$, so that $k'_{j}=k_{j}-r_{j}$.
Thus, noting that all $|\lambda_{w}|\leq1$,
\begin{align*}
\left|\langle S_{f}(g),DS_{f}(h)\rangle\right| & \le\sum_{(k_{0},\dots,k_{p})\in K}\prod_{j=0}^{p}|c_{j}(k_{j})c_{j}(k_{j}-r_{j})|\\
 & =\prod_{j=0,j\neq m}^{p}\sum_{k_{j}\in\mathbb{Z}}|c_{j}(k_{j})c_{j}(k_{j}-r_{j})|\cdot\sum_{k_{m}\in\mathbb{Z}}|c_{m}(k_{m})c_{m}(k_{m}-r_{m})|\\
 & \leq\sum_{k}|c_{m}(k)c_{m}(k-r_{m})|,
\end{align*}
since by Cauchy-Schwarz, $\sum_{k_{j}}|c_{j}(k_{j})c_{j}(k_{j}-r_{j})|\leq\Vert\{|c_{j}(k_{j})|\}_{k_{j}}\Vert_{\ell^{2}}\Vert\{|c_{j}(k_{j}-r_{j})|\}_{k_{j}}\Vert_{\ell^{2}}\leq\Vert f^{d_{j}}\Vert_{2}^{2}=1$. 

Since $r_{m}\neq0$, $d_{m}\in\{\pm1\}$; thus the coefficients $c_{m}(k)$
are either coefficients of $f$ or $\bar{f}$. Thus, denoting by $a$
either $r_{m}$ or $-r_{m}$, 
\[
\left|\langle S_{f}(g),DS_{f}(h)\rangle\right|\leq\sum_{k}|c(k)c(k-a)|,
\]
where $c(k)=\int f(z)z^{-k}d\mu(z)$ are the Fourier coefficients
of $f$. Then $\eta_{f}=\sum_{k\neq0}|c(k)|^{2}$. In particular,
$|c(k)|\leq\sqrt{\eta_{f}}$ if $k\neq0$ (and for any $k$, $|c(k)|\leq1$).
Thus we can write
\begin{align*}
\sum_{k}|c(k)c(k-a)| & =|c(0)c(-a)|+|c(a)c(0)|+\sum_{k\neq0,k\neq a}|c(k)c(k-a)|\\
 & \leq2\sqrt{\eta_{f}}+\left(\sum_{k\neq0,k\neq a}|c(k)|^{2}\right)^{1/2}\left(\sum_{k\neq0,k\neq a}|c(k-a)|^{2}\right)^{1/2}\leq2\sqrt{\eta_{f}}+\eta_{f},
\end{align*}
as claimed. 
\end{proof}
The following Lemma describes a set of automorphisms that will be
used to satisfy the hypothesis of Proposition \ref{prop:TandTheta}.
\begin{lem}
\label{lem:Psifm}Let $f:\mathbb{T\to\mathbb{T}}$ be as before. Let
\[
\Psi_{f,m}=R_{-m}\circ S_{f}\circ R_{m}.
\]
Then $\Psi_{f,m}$ is a normal trace-preserving automorphism of $M$,
$\Psi_{f,m}(w)=w$, $\Psi_{f,m}(A)=Af(BA^{-m})$ and $\Vert\Psi_{f,m}(A)-A\Vert_{2}=\Vert f-1\Vert_{L^{2}(\mathbb{T})}.$
Moreover, $\Psi_{f,m}^{-1}=\Psi_{\bar{f},m}$.
\end{lem}
\begin{proof}
Since both $S_{f}$ and $R_{m}$ preserve $w$, so does $\Psi_{f,m}$.
Furthermore, $S_{f}(R_{m}(A))=S_{f}(A)=Af(B)$. Thus $R_{-m}(S_{f}(R_{m}(A)))=R_{-m}(Af(B))=Af(BA^{-m})$.
Finally, $\Vert\Psi_{f,m}(A)-A\Vert_{2}=\Vert Af(BA^{-m})-A\Vert_{2}=\Vert f(BA^{-m})-1\Vert_{2}=\Vert f-1\Vert_{L^{2}(\mathbb{T})}$
since $BA^{-m}$ is a Haar unitary.
\end{proof}

\subsection{Bass-Serre decomposition of $\mathbb{F}_{2}$ . }

Let $H$ denote the subgroup generated by $A$ and $w$. As before,
set $N=W^{*}(A,w)=L(H)$. Note that the automorphism $R_{m}$ satisfies
$R_{m}(A)=A$ and $R_{m}(w)=w$, so $R_{m}(h)=h$ for any $h\in H$. 

For an element $g\in\mathbb{F}_{2}$ let $\ell(g)$ be the smallest
number of letters $B$ or $B^{-1}$ required to write $g$ as a word
in $B,B^{-1}$ and arbitrary elements of $H$. In particular, $\ell(g)=0$
iff $g\in H$. 

Let $C=BAB^{-1}$. Then $H$ is the subgroup generated by $A$ and
$C$. Moreover, $\mathbb{F}_{2}$ is the HNN-extension of $H$: it
can be viewed as freely generated by letters $A,C,B$ subject to the
relation $BAB^{-1}=C$. The associated Bass-Serre tree $\mathcal{T}$
has as vertex set $\mathbb{F}_{2}/H$, and edges between $xH$ and
$xBH$ for any $x\in\mathbb{F}_{2}$. Thus paths in the tree correspond
to words in $B,B^{-1}$ and arbitrary elements of $H$. Thus the geodesic
distance $d_{\mathcal{T}}$ on $\mathcal{T}$ satisfies
\[
d_{\mathcal{T}}(xH,yH)=\ell(x^{-1}y).
\]

For two probability measures $\mu_{1},\mu_{2}$ on $\mathbb{F}_{2}$
define their convolution $\mu_{1}*\mu_{2}$ by
\[
\mu_{1}*\mu_{2}(g)=\sum_{g=g_{1}g_{2}}\mu_{1}(g_{1})\mu_{2}(g_{2})
\]
and total variation distance
\[
\Vert\mu_{1}-\mu_{2}\Vert_{TV}=\sup_{E\subset\mathbb{F}_{2}}|\mu_{1}(E)-\mu_{2}(E)|=\frac{1}{2}\sum_{g\in\mathbb{F}_{2}}|\mu_{1}(g)-\mu_{2}(g)|.
\]
Note that if $\mu_{1}$ and $\mu_{2}$ are probability measures (thus
with total mass $1$), then
\begin{align*}
\Vert\mu_{1}-\mu_{2}\Vert_{TV} & =\frac{1}{2}\sum_{g\in\mathbb{F}_{2}}|\mu_{1}(g)-\mu_{2}(g)|\\
 & =\frac{1}{2}\sum_{g\in\mathbb{F}_{2}}(\mu_{1}(g)+\mu_{2}(g)-2\min(\mu_{1}(g),\mu_{2}(g)))\\
 & =1-\sum_{g\in\mathbb{F}_{2}}\min(\mu_{1}(g),\mu_{2}(g)).
\end{align*}

Given $g\in\mathbb{F}_{2}$, $f:\mathbb{T}\to\mathbb{T}$ measurable,
and $m\in\mathbb{Z}$, define the probability measure $\hat{\mu}_{g}^{f,m}$
on $\mathbb{F}_{2}$ by
\[
\hat{\mu}_{g}^{f,m}(h)=\left|\langle h,\Psi_{f,m}(g)\rangle\right|^{2}.
\]

The following lemma shows that if $m$ is large enough, then the measures
$\hat{\mu}_{g}^{f.m}$ are, up to a factor that only depends on $g$
and on how close if $f$ to $1$, given by a convolution.
\begin{lem}
\label{lem:TotalVariationalEstimate}Let $g\in\mathbb{F}_{2}$ and
choose a shortest relative expression for $g$ corresponding to a
geodesic from $gH$ to the identity element $H$ in $\mathcal{T}$:
\[
g=h_{0}B^{\epsilon_{1}}h_{1}B^{\epsilon_{2}}\cdots B^{\epsilon_{l}}h_{l},\qquad\epsilon_{j}\in\{\pm1\},\ h_{j}\in H,\ l=\ell(g).
\]
Let
\[
\hat{\nu}_{g}^{f,m}=\delta_{h_{0}}*\hat{\mu}_{B^{\epsilon_{1}}}^{f,m}*\delta_{h_{1}}*\cdots*\hat{\mu}_{B^{\epsilon_{l}}}^{f,m}*\delta_{h_{l}}
\]
(with the convention that $\hat{\nu}_{g}^{f,m}=\delta_{g}$ if $g\in H$).
Then there are constants $C_{g}$ and $m_{g}$ so that for any $m>m_{g}$
and every measurable $f:\mathbb{T}\to\mathbb{T}$, 
\[
\Vert\hat{\nu}_{g}^{f,m}-\hat{\mu}_{g}^{f,m}\Vert_{TV}\leq C_{g}\eta_{f},
\]
with $\eta_{f}=1-|\tau(f(B))|^{2}$. 
\end{lem}
\begin{proof}
Recall that $\Psi_{f,m}=R_{-m}\circ S_{f}\circ R_{m}$. Let
\[
\mu_{g}^{f,m}(h)=\left|\langle h,S_{f}R_{m}(g)\rangle\right|^{2}
\]
and
\[
\nu_{g}^{f,m}=\delta_{h_{0}}*\mu_{B^{\epsilon_{1}}}^{f,m}*\delta_{h_{1}}*\cdots*\mu_{B^{\epsilon_{l}}}^{f,m}*\delta_{h_{l}}
\]

Since $R_{-m}$ is a group automorphism that fixes $H$ pointwise,
it preserves convolutions, total variation distance, and measures
$\delta_{h_{j}}$, and thus
\[
\mu_{g}^{f,m}=(R_{m})_{*}\hat{\mu}_{g}^{f,m},\qquad\nu_{g}^{f,m}=(R_{m})^{*}\hat{\nu}_{g}^{f,m}.
\]
Thus the claimed inequality is equivalent to proving that for some
$C_{g}$ and $m_{g}$, and arbitrary $f$ and $m>m_{g}$,
\[
\Vert\nu_{g}^{f,m}-\mu_{g}^{f,m}\Vert_{TV}\leq C_{g}\eta_{f}.
\]

We first consider the case that $l=\ell(g)>0$.

Recall that by assumption $g$ has a minimal representation
\[
g=h_{0}B^{\epsilon_{1}}h_{1}B^{\epsilon_{2}}\cdots B^{\epsilon_{l}}h_{l},\qquad\epsilon_{j}\in\{\pm1\},\ h_{j}\in H,\ l=\ell(g).
\]
If $\epsilon_{j}\neq\epsilon_{j+1}$, then minimality implies that
$h_{j}$ is not a power of $A$ (if $\epsilon_{j}=1$) or not a power
of $BAB^{-1}$ (otherwise). 

Since $R_{m}$ fixes $H$ pointless, $R_{m}(g)$ is obtained form
$g$ by replacing each $B$ by $BA^{m}$ and $B^{-1}$ by $A^{-m}B^{-1}$
. The resulting word need not be reduced: it has the form
\[
k_{0}A^{\epsilon_{1}m}k_{1}A^{\epsilon_{2}m}\cdots A^{\epsilon_{l}m}k_{l}
\]
where for $0<j<l$ the words $k_{j}$ are given by:
\[
k_{j}=\begin{cases}
h_{j}B, & \epsilon_{j}=\epsilon_{j+1}=1\\
h_{j}, & \epsilon_{j}=1,\ \epsilon_{j+1}=-1\\
B^{-1}h_{j}B, & \epsilon_{j}=-1,\ \epsilon_{j+1}=1\\
B^{-1}h_{j} & \text{otherwise.}
\end{cases}
\]
and
\begin{align*}
k_{0} & =\begin{cases}
h_{0}B, & \epsilon_{1}=1\\
h_{0} & \textrm{otherwise;}
\end{cases} & k_{l}=\begin{cases}
h_{l}, & \epsilon_{l}=1\\
B^{-1}h_{l} & \textrm{otherwise.}
\end{cases}
\end{align*}
Thus cancellation is possible between the start of $k_{j}$ or the
end of $k_{j-1}$ and $A^{\epsilon_{j}m}$ (note that if $0<j<l$,
$k_{j}$ must contain at least one letter $B$, i.e., it has the form
$k_{j}=A^{a_{j}}v_{j}A^{b_{j}}$ with $v_{j}$ nonempty and starting
and ending with a power of $B$; and $k_{0}=v_{0}A^{b_{0}}$, $k_{l}=A^{a_{l}}v_{l}$,
where $v_{0}$ may be empty or must end with $B^{\pm1}$, and $v_{l}$
must be empty or must start with $B^{\pm1}$). If we set
\[
K=1+2l+\sum_{j=0}^{l}|h_{j}|
\]
(where we put $|\cdot|$ for the word metric on $\mathbb{F}_{2}$
associated to the generating set $A,B$), then there cannot be more
than $K$ cancellations between $k_{j}$ or $k_{j-1}$ and $A^{\epsilon_{j}m}$.
Thus if we choose $m>2K+1$, we can ensure that after reduction,
\[
R_{m}(g)=v_{0}A^{t_{1}}v_{1}\cdots A^{t_{l}}v_{l}
\]
where
\[
t_{j}=b_{j-1}+\epsilon_{j}m+a_{j}.
\]
and same $v_{0},\dots,v_{l}$ as before. Note that the sign of $t_{j}$
is $\epsilon_{j}$. 

Write each block $A^{t_{j}}$ as $A^{\epsilon_{j}K+b_{j-1}}\cdot A^{\epsilon_{j}(m-2K)}\cdot A^{\epsilon_{j}K+a_{j}}=S_{j}M_{j}E_{j}$.
Then
\[
R_{m}(g)=v_{0}S_{1}M_{1}E_{1}v_{1}S_{2}M_{2}E_{2}\cdots v_{l}
\]
is a reduced word. 

Let $r=m-2K$. For each multi-index $\mathbf{k}_{j}=(k_{j,1},\dots,k_{j,r-1})$
put
\[
M_{j}(\mathbf{k}_{j})=A^{\epsilon_{j}}B^{k_{j,1}}A^{\epsilon_{j}}B^{k_{j,2}}\cdots B^{k_{j,r-1}}A^{\epsilon_{j}}
\]
 and define
\[
V(\mathbf{k})=v_{0}S_{1}M_{1}(\mathbf{k}_{1})E_{1}\cdots v_{l}.
\]
Define also
\[
U_{j}(\mathbf{k}_{j})=\begin{cases}
BA^{K}M_{j}(\mathbf{k}_{j})A^{K} & \epsilon_{j}=1\\
A^{-K}M_{j}(\mathbf{k}_{j})A^{-K}B^{-1} & \textrm{otherwise}.
\end{cases}
\]
Then for any $\mathbf{k}$,
\[
V(\mathbf{k})=h_{0}U_{1}(\mathbf{k}_{1})h_{1}\cdots U_{l}(\mathbf{k}_{l})h_{l}.
\]
Indeed, in the process of reducing $h_{0}U_{1}(\mathbf{k}_{1})h_{1}\cdots U_{l}(\mathbf{k}_{l})h_{l}$
we can perform the same sequence of cancellations as when reducing
$\ensuremath{h_{0}R_{m}(B^{\epsilon_{1}})h_{1}\cdots R_{m}(B^{\epsilon_{l}})h_{l}}$,
and none of these cancellations occur inside the words $M_{j}(\mathbf{k}_{j})$. 

For each $\mathbf{k}$ let
\[
F(\mathbf{k})=\prod_{j=1}^{l}\prod_{t=1}^{r-1}|c_{\epsilon_{j}}(k_{j,t})|^{2}
\]
where as before $c_{\epsilon}(k)$ is the $k$-th Fourier coefficient
of $f^{\epsilon}$. 

Let $p$ be the total number of occurrences of $A^{\pm1}$ in $V(\mathbf{k})$,
and let $d_{0},\dots,d_{p}$ be as described in the paragraph before
Lemma \ref{lem:descriptionOfSfWords}. Note that if the $i$-th occurrence
of $A^{\pm1}$ happens inside $M_{j}$ and not at the end of $M_{j}$,
then $d_{i}=\epsilon_{j}\neq0$. Thus if we denote by $N$ the total
number of nonzero $d_{j}$, then $N=l(r-1)+N_{0}$. 

If we put $S=\sum_{j=0}^{l}|h_{j}|,$then $p\leq lm+S$ and therefore
$N\leq lm+S+1$. Thus 
\[
N_{0}=N-l(r-1)\leq lm+S+1-l(r-1)=S+1+(2K+1)l.
\]

Lemma \ref{lem:descriptionOfSfWords} then implies that
\[
\mu_{g}^{f,m}(V(\mathbf{k}))=F(\mathbf{k})|c(0)|^{2N_{0}}
\]

Similarly, 
\[
\mu_{B^{\epsilon_{j}}}^{f,m}(U_{j}(\mathbf{k}_{j}))=\left(\prod_{t=1}^{r-1}|c_{\epsilon_{j}}(k_{j,t})|^{2}\right)\cdot|c(0)|^{2(2K+1)}\geq0
\]
and so by definition of convolution and since $V(\mathbf{k})=h_{0}U_{1}(\mathbf{k}_{1})h_{1}\cdots U_{l}(\mathbf{k}_{l})h_{l}$,
$\nu_{g}^{f,m}(V(\mathbf{k}))$ is a sum of non-negative terms, at
least one of which is 
\[
\prod_{j=1}^{l}\mu_{B^{\epsilon_{j}}}^{f,m}(U_{j}(\mathbf{k}_{j}))=F(\mathbf{k})|c(0)|^{2N_{0}'},\qquad N_{0}'=(2K+1)l.
\]
Thus
\[
\nu_{g}^{f,m}(V(\mathbf{k}))\geq F(\mathbf{k})|c(0)|^{2N_{0}'}.
\]

In particular, for any $\mathbf{k}$,
\[
\min(\nu_{g}^{f,m}(V(\mathbf{k})),\mu_{g}^{f,m}(V(\mathbf{k})))\geq|c(0)|^{2\max(N_{0},N_{0}')}F(\mathbf{k}).
\]

Since $N_{0}\leq S+1+(2K+1)l$ and $N_{0}'=(2K+1)l$, $\max(N_{0},N_{0}')\leq S+1+(2K+1)l$. 

Let $c(0)=\hat{f}(0)$. Because the map $\mathbf{k}\to V(\mathbf{k})$
is injective,
\begin{align*}
\sum_{h\in\mathbb{F}_{2}}\min(\nu_{g}^{f,m}(h),\mu_{g}^{f,m}(h)) & \geq\sum_{\mathbf{k}}\min(\nu_{g}^{f,m}(V(\mathbf{k})),\mu_{g}^{f,m}(V(\mathbf{k})))\\
 & \geq\sum_{\mathbf{k}}F(k)|c(0)|^{2\max(N_{0},N_{0}')}\\
 & =|c(0)|^{2\max(N_{0},N_{0}')}\prod_{j=1}^{l}\prod_{t=1}^{r-1}\sum_{k_{j,t}}|c_{\epsilon_{j}}(k_{j,t})|^{2}\\
 & =|c(0)|^{2\max(N_{0},N_{0}')}\prod_{j=1}^{l}\prod_{t=1}^{r-1}\Vert f^{\epsilon_{j}}\Vert_{L^{2}}^{2}=|c(0)|^{2\max(N_{0},N_{0}')}.
\end{align*}
\begin{align*}
\Vert\hat{\nu}_{g}^{f,m}-\hat{\mu}_{g}^{f,m}\Vert_{TV} & \leq1-|c(0)|^{2\max(N_{0},N_{o}')}\\
 & \leq(1+|c(0)|^{2}+\cdots+|c(0)|^{2\max(N_{0},N_{0}')-2})(1-|c(0)|^{2})\\
 & \leq(S+1+(2K+1)l)\eta_{f}=C_{g}\eta_{f},
\end{align*}
with the constant $C_{g}$ depending only on $g$. 

Let us now consider the case $\ell(g)=0$, i.e., $g\in H$. Then $\nu_{g}^{f,m}=\delta_{g}$.
Since $g\in H$, $g$ can be written a s product of $C=BAB^{-1}$
and $A$. Thus
\[
g=B^{n_{0}}A^{\epsilon_{1}}B^{n_{1}}A^{\epsilon_{2}}\cdots A^{\epsilon_{p}}B^{n_{p}}
\]
In this case \ref{lem:descriptionOfSfWords} shows that $\mu_{g}^{f,m}(g)\geq|c(0)|^{2r_{g}}$
where $r_{g}$ is the total number of letters $A$ in the expansion
above. Thus $\Vert\mu_{g}^{f,m}-\delta_{g}\Vert_{TV}\leq1-|c(0)|^{2r_{g}}\leq r_{g}\eta_{f}$.
\end{proof}

\subsection{The Operators $T_{q}$ and Quadratic Forms $Q_{q}$ .}

For each $0\leq q\leq1$, let
\[
T_{q}:L^{2}(M)\to L^{2}(M)
\]
be the diagonal operator
\[
T_{q}\left(\sum\alpha_{g}g\right)=\sum\alpha_{g}q^{\ell(g)}g,
\]
Then $T_{q}$ is a positive contraction. For $q=0$, $T_{0}=E_{N}$,
and as $q\to1$, $T_{q}\to\operatorname{Id}$ strongly. The associated
quadratic form $Q_{k}(x)=\langle T_{q}x,x\rangle$ is given by 
\[
Q_{q}(\sum\alpha_{g}g)=\sum|\alpha_{g}|^{2}q^{\ell(g)}.
\]
with the convention that $0^{0}=1$. 

Note that since $R_{m}$ fixes $H$ pointless and replaces each $B^{\pm1}$
by some word containing only one $B^{\pm1}$, $\ell(R_{m}(g))\le\ell(g)$.
Applying the same logic to $R_{-m}=R_{m}^{-1}$ shows that $\ell(R_{m}(g))=\ell(g)$.
In particular, this means that $T_{q}$ commutes with the unitary
implementing $R_{m}$ on $L^{2}(M)$. 
\begin{cor}
\label{cor:lowerEstimatePsigVsPsiB}For any $g\in\mathbb{F}_{2}$,
there are constants $C_{g},m_{g}$ so that whenever $f:\mathbb{T}\to\mathbb{T}$
is measurable, $\eta_{f}=1-|\tau(f(B))|^{2}$, for any $m>m_{g}$
and any $0\leq q\leq1$,
\[
Q_{q}(\Psi_{f,m}^{-1}(g))\geq Q_{q}(\Psi_{f,m}^{-1}(B))^{\ell(g)}-C_{g}\eta_{f}.
\]
\end{cor}
\begin{proof}
By definition,
\[
Q_{q}(\Psi_{f,m}^{-1}(g))=\sum_{h}q^{\ell(h)}\hat{\mu}_{g}^{\bar{f},m}=\sum_{h}q^{\ell(h)}\mu_{g}^{\bar{f},m}(h).
\]
Since $q\in[0,1]$, the Total Variation estimate from Lemma \ref{lem:TotalVariationalEstimate}
implies that
\begin{align*}
Q_{q}(\Psi_{f,m}^{-1}(g)) & \geq\sum_{h}q^{\ell(h)}\hat{\nu}_{g}^{\bar{f},m}(h)-C_{g}\eta_{f}\geq\prod_{j=1}^{\ell(g)}Q_{q}(\Psi_{f,m}^{-1}(B^{\epsilon_{j}}))-C_{g}\eta_{f}\\
 & \geq Q_{q}(\Psi_{f,m}^{-1}(B))^{\ell}-C_{g}\eta_{f},
\end{align*}
where the last inequality comes from the definition of $\nu_{g}^{f,m}(h)$
and the observation that if $h_{j}\in H$, then $\ell(h_{0}s_{1}h_{1}\dots s_{\ell(g)}h_{\ell(g)})\leq\sum_{j=1}^{\ell(g)}\ell(s_{j})$
\end{proof}

\subsection{Choosing $f$ and $m$. }

For $t>0$ and $m\geq1$ let
\[
f_{m,t}(z)=\exp\left[i\sqrt{\frac{t}{2m}}(z+z^{-1})\right],\qquad\theta_{m,t}=\Psi_{f_{m,t},m}\in\operatorname{Aut}(M).
\]
Since $|\exp(iz)-1|\leq|z|$, we have that
\[
|f_{m,t}(z)-1|\leq\sqrt{\frac{t}{m}}\cdot\left|\frac{z+z^{-1}}{\sqrt{2}}\right|.
\]
Thus
\begin{align*}
\Vert f_{m,t}-1\Vert_{L^{2}(\mathbb{T})} & =\left[\int|f_{m,t}(z)-1|^{2}d\mu(z)\right]^{1/2}\\
 & \leq\sqrt{t/m}\left[\int\left|\frac{z+z^{-1}}{\sqrt{2}}\right|^{2}d\mu(z)\right]^{1/2}=\sqrt{t/m}\left[\int2\cos^{2}2\pi(t)dt\right]^{1/2}=\sqrt{t/m},
\end{align*}
In particular,
\[
\eta_{f_{m,t}}\leq\Vert f_{m,t}-1\Vert^{2}\leq t/m.
\]

We have made this concrete choice of $f_{m,t}$ for definiteness.
In fact, all we need are the following properties of $f_{m,t}$:
\begin{lem}
\label{lem:propertiesOffmt}(i) $f_{m,t}:\mathbb{T}\to\mathbb{T}$
is measurable. (ii) Let $\ensuremath{c_{m}(k)=\int_{\mathbb{T}}\bar{f}_{m,t}(z)z^{-k}d\mu(z)=\widehat{\bar{f}_{m,t}}(k)}$,
and let $\alpha_{m}=\left|\widehat{\bar{f}_{m,t}}(0)\right|^{2}$,
$\beta_{m}=\left|\widehat{\bar{f}_{m,t}}(\pm1)\right|^{2}$. Then
$\alpha_{m}=1-(t/m)+O(m^{-2})$ and $\beta_{m}=(t/2m)+O(m^{-2})$.
(iii) The Fourier series for $f$ is absolutely convergent. 
\end{lem}
\begin{proof}
Taylor expansion gives $\alpha_{m}=1-(t/m)+O(m^{-2})$ and $\beta_{m}=(t/2m)+O(m^{-2})$. 
\end{proof}
\begin{lem}
The automorphism $\theta_{m,t}$ satisfies: $\theta_{m,t}(w)=w$ and
$\Vert\theta_{m,t}(A)-A\Vert_{2}\leq\sqrt{t/m}.$
\end{lem}
\begin{proof}
Both follow immediately from Lemma \ref{lem:Psifm}. 
\end{proof}
For $t>0$ and $q\in[0,1]$ let 
\[
F_{t}(q)=e^{t(q-1)/2}\left(\frac{1-e^{-t}}{2}+\frac{1+e^{-t}}{2}q\right).
\]
The function $F_{t}(q)$ will be later shown to be the limit of certain
generating functions used to describe the combinatorics of $\theta_{m,t}^{-1}(B)$. 
\begin{lem}
\label{lem:F_tlowerEqThetaOfB}For every $t>0$ and every $\epsilon>0$
there exists an $m_{0}$ so that for all $m>m_{0}$ and all $0\leq q\leq1$,
\[
Q_{q}(\theta_{m,t}^{-1}(B))\geq F_{t}(q)-\epsilon.
\]
\end{lem}
\begin{proof}
Let $\ensuremath{c_{m}(k)=\int_{\mathbb{T}}\bar{f}_{m,t}(z)z^{-k}d\mu(z)=\widehat{\bar{f}_{m,t}}(k)}$,
and let $\alpha_{m}=\left|\widehat{\bar{f}_{m,t}}(0)\right|^{2}$,
$\beta_{m}=\left|\widehat{\bar{f}_{m,t}}(\pm1)\right|^{2}$ as in
Lemma \ref{lem:propertiesOffmt}; thus $\alpha_{m}=1-(t/m)+O(m^{-2})$
and $\beta_{m}=(t/2m)+O(m^{-2})$. 

Since $\theta_{m,t}^{-1}=R_{-m}\circ S_{\bar{f}_{m,t}}\circ R_{m}$,
it follows from commutation between $T_{q}$ and $R_{m}$ that $Q_{q}(\theta_{m,t}^{-1}(B))=Q_{q}(S_{\bar{f}_{m,t}}(R_{m}(B)))=Q_{q}(S_{\bar{f}_{m,t}}(BA^{m}))$.

Lemma \ref{lem:descriptionOfSfWords} then gives us
\[
S_{\bar{f}_{m,t}}(BA^{m})=\sum_{\mathbf{n}\in\mathbb{Z}^{m}}\left(\prod_{j=1}^{m}c_{m}(n_{j})\right)W(\mathbf{n})
\]
where $\mathbf{n}\to W(\mathbf{n})$ is the injective map
\[
W(\mathbf{n})=BAB^{n_{1}}AB^{n_{2}}\cdots AB^{n_{m}}.
\]
Thus
\begin{align}
Q_{q}(\theta_{m,t}^{-1}(B)) & =\sum_{\mathbf{n}\in\mathbb{Z}^{m}}\left(\prod_{j=1}^{m}|c_{m}(n_{j})|^{2}\right)q^{\ell(W(\mathbf{n}))}\nonumber \\
 & \geq\sum_{\mathbf{n}\in\{-1,0,1\}^{m-1}\times\{0\}}\left(\prod_{j=1}^{m}|c_{m}(n_{j})|^{2}\right)q^{\ell(W(\mathbf{n}))}.\label{eq:termsWeKeep}
\end{align}
We will use the smaller sum (keeping only terms with $n_{1},\dots,n_{m-1}\in\{-1,0,1\}$
and $n_{m}=0$) for our lower estimate. Thus from now on we only consider
$n_{m}=0$ and each $n_{1},\dots,n_{m-1}\in\{-1,0,1\}$. Our aim is
to bound $\ell(W(\mathbf{n})).$

Recall that $H$ is generated by $A$ and $C=BAB^{-1}$. 

Define $P_{0}=BA$, $P_{j}=P_{j-1}B^{n_{j}}A$, $1\leq j\leq m-1$.
Since $n_{m}=0$,
\[
W(\mathbf{n})=P_{m-1}.
\]
We recursively factor
\[
P_{j}=V_{j}B^{e_{j}},\qquad e_{j}\in\{0,1\}
\]
as follows. Initially, $P_{0}=BA=CB$, thus we set $V_{0}=C$ and
$e_{0}=1$. Suppose we factored $P_{j}=V_{j}B^{e_{j}}$. Then
\[
P_{j+1}=P_{j}B^{n_{j+1}}A.
\]
We then set $V_{j+1},e_{j+1}$ according to the table:
\[
(V_{j+1},e_{j+1})=\begin{cases}
(V_{j}B^{-1}A,0), & e_{j}+n_{j+1}=-1,\\
(V_{j}A,0) & e_{j}+n_{j+1}=0,\\
(V_{j}C,1) & e_{j}+n_{j+1}=1,\\
(V_{j}BC,1) & e_{j}+n_{j+1}=2.
\end{cases}
\]
Because $BA=CB$ and $B^{2}A=BCB$, $P_{j+1}=V_{j+1}B^{e_{j+1}}$.
We claim that then $P_{j+1}=V_{j+1}B^{e_{j+1}}$. Indeed, there are
6 possible steps:
\begin{itemize}
\item (Type 1) Suppose that $e_{j}+n_{j+1}=-1$ , so that $e_{j}=0$ and
$n_{j+1}=-1$. Then $P_{j}=V_{j}B^{0}$ , $V_{j+1}=V_{j}B^{-1}A$,
and $P_{j+1}=P_{j}B^{-1}A=V_{j+1}B^{0}$.
\item Suppose that $e_{j}+n_{j+1}=0$. There are two possibilities:
\begin{itemize}
\item (Type 2) If $e_{j}=0$ and $n_{j+1}=0$ , then $P_{j}=V_{j}B^{o}$,
$V_{j+1}=V_{j}A$, and $P_{j+1}=P_{j}A=V_{j+1}B^{0}.$
\item (Type 3) If $e_{j}=1$ and $n_{j+1}=-1$, then $P_{j}=V_{j}B$, $V_{j+1}=V_{j}A$,
and $P_{j+1}=P_{j}B^{-1}A=V_{j}A=V_{j+1}B^{0}.$
\end{itemize}
\item Suppose that $e_{j}+n_{j+1}=1$. There are two possibilities:
\begin{itemize}
\item (Type 4) If $e_{j}=0$ and $n_{j+1}=1$, then $P_{j}=V_{j}B^{0}$,
$V_{j+1}=V_{j}C=V_{j}BAB^{-1}$, and $P_{j+1}=P_{j}BA=V_{j+1}B.$
\item (Type 5) If $e_{j}=1$ and $n_{j+1}=0$, then $P_{j}=V_{j}B$, $V_{j+1}=V_{j}C=V_{j}BAB^{-1}$,
and $P_{j+1}=P_{j}A=V_{j+1}B.$
\end{itemize}
\item (Type 6) Suppose finally that $e_{j}+n_{j+1}=2$ , so that $e_{j}=1$
and $n_{j+1}=1$ . Then $P_{j}=V_{j}B$, $V_{j+1}=V_{j}BC=V_{j}B^{2}AB^{-1},$and
$P_{j+1}=P_{j}BA=V_{j+1}B.$
\end{itemize}
Let us say that the factorization $P_{j}=V_{j}B^{e_{j}}$ is even
(respectively, odd) if $e_{j}=0$ (resp., $e_{j}=1$).

If the transition from $P_{j}=V_{j}B^{e_{j}}$ to $P_{j+1}=V_{j+1}B^{e_{j+1}}$
is of type 2 or 5, then the parity is unchanged, and $\ell(V_{j+1})=\ell(V_{j})$
since $V_{j}$ is multiplied by an element of $H$. We call these
transitions \textbf{type A}. In this case, $n_{j+1}=0$.

If the transition is of type 3 or 4, then the parity is reversed,
but $\ell(V_{j+1})=\ell(V_{j})$ since $V_{j}$ is multiplied by an
element of $H$. We call these transitions \textbf{type B}. In this
case, $|n_{j+1}|=1$.

Finally, if the transition is of type 1 or 6, then the parity is unchanged,
but $V_{j}$ is changed by the multiplication with either $B$ or
$B^{-1}$ as well as an element of $H$. Thus the length of $V_{j+1}$
could have increased by $1$. We call these transitions \textbf{type
C}. In this case, $|n_{j+1}|=1$.

Let $a_{j}$ , $b_{j}$, and $c_{j}$ be the total number of type
A, type B, and type C transitions that occurred when moving from $P_{0}=V_{0}B$
to $P_{j}=V_{j}B^{e_{j}}$ . Then
\[
\ell(V_{j})\leq c_{j},\qquad\ell(P_{j})\leq c_{j}+e_{j}.
\]
Let $a=a_{m-1}$, $b=b_{m-1}$, $c=c_{m-1}$ be the total number of
transitions of each type. Then $a+b+c=m-1$. Since only type B transitions
affect parity, and the inital parity is odd, the final parity is odd
iff $b$ is even. 
\[
\ell(W(\mathbf{n}))=\ell(P_{m-1})\leq c+1_{b\textrm{ even}}.
\]
Moreover, the number of indices $j$ for which $|n_{j}|=1$ is precisely
$c+b$ and the number of $j$ for which $n_{j}=0$ is precisely $a$.
Therefore, the term in (\ref{eq:termsWeKeep}) associated to this
particular $\mathbf{n}$ is 
\begin{equation}
\alpha_{m}^{a+1}\beta_{m}^{b+c}q^{\ell(W(\mathbf{n}))}\geq\alpha_{m}^{a+1}\beta_{m}^{b+c}q^{c+1_{b\textrm{ even}}}.\label{eq:estimateOnabTerm}
\end{equation}

Note also that if we are given the sequence of transition types, then
we can recover $\mathbf{n}$ fully. Indeed, the value of $n_{j+1}$
is completely determined if we know the transition type and the parity
$e_{j}$: for transitions of type A, $n_{j+1}=0$; for type B, $n_{j+1}$
is negative iff $e_{j}=1$, and for type C, $n_{j+1}$ is positive
iff $e_{j}=1$. Since $e_{0}=1$, it follows that the sequence of
transition types completely determines $\mathbf{n}$. Conversely,
knowing $\mathbf{n}$ produces a sequence of $m-1$ transition types.
Therefore, there is a bijection between sequences of $m-1$ transition
types and tuples $\mathbf{n}$. One can realize this transition explicitly
using a generating function.

Consider the polynomial $\mathscr{P}_{m}(X,Y)$ given by
\[
\mathscr{P}_{m}(X,Y)=(\alpha_{m}+\beta_{m}X+\beta_{m}Y)^{m-1}\alpha_{m}.
\]
If one expands the right hand side of this equation, then each term
is a product of $m$ terms, each chosen from among $\alpha_{m}$,
$\beta_{m}X$ and $\beta_{m}Y$. Keeping the order of terms in this
product, we see that the products are in bijection with sequences
of $m-1$ of three types, with the three terms corresponding to types
A, B, and C. Moreover, $b$ and $c$ correspond precisely to the total
power of $X$ and $Y$ in the corresponding term.

Let $p_{b,c}^{(m)}$ be the coefficient of $X^{b}Y^{c}$ in $\mathscr{P}_{m}(X,Y)$.
Then
\[
p_{b,c}^{(m)}=\sum_{\mathbf{n}\textrm{ has \ensuremath{b} type B and \ensuremath{c} type C}}\alpha_{m}^{a+1}\beta_{m}^{b+c},
\]
where $a=m-b-c-1.$ Therefore, (\ref{eq:termsWeKeep}) combined with
(\ref{eq:estimateOnabTerm}) implies that
\[
Q_{q}(\theta_{m,t}^{-1}(B))\geq q\sum_{b,c\geq0\textrm{ with \ensuremath{b} even}}p_{b,c}^{(m)}q^{c}+\sum_{b,c\geq0\textrm{ with \ensuremath{b} odd}}p_{b,c}^{(m)}q^{c}.
\]
These sums can be recovered by considering the even and odd parts
of $\mathscr{P}$ and substituting $X=q$ and $Y=\pm1$: 
\[
\sum_{b,c\textrm{ with \ensuremath{b} even}}p_{b,c}^{(m)}q^{c}=\frac{1}{2}(\mathscr{P}_{m}(1,q)+\mathscr{P}_{m}(-1,q))
\]
and
\[
\sum_{b,c\textrm{ with \ensuremath{b} odd}}p_{b,c}^{(m)}q^{c}=\frac{1}{2}(\mathscr{P}_{m}(1,q)-\mathscr{P}_{m}(-1.-q)).
\]
Substituting this in gives us the final estimate
\[
Q_{q}(\theta_{m,t}^{-1}(B))\geq\frac{q+1}{2}\mathscr{P}_{m}(1,q)+\frac{q-1}{2}\mathscr{P}_{m}(-1,q).
\]

Let $\sigma\in\{-1,1\}$. We now substitute $\alpha_{m}=1-(t/m)+O(m^{-2})$
and $\beta_{m}=(t/2m)+O(m^{-2})$ into the equation for $\mathscr{P}_{m}$:
\[
\mathscr{P}_{m}(q,\sigma)=(\alpha_{m}+\beta_{m}q+\beta_{m}\sigma)^{m-1}\alpha_{m}.
\]
For $m$ large enough $\alpha_{m}+\beta_{m}(q+\sigma)=1-(t/m)+(t/2m)(q+\sigma)+O(m^{-2})$
(where $O(m^{-2})$ is uniform in $q$). In particular, for large
enough $m$ this is positive for all $q\in[-1,1]$ and both choices
of $\sigma$. Taking logarithms of both sides gives us
\begin{align*}
\log\mathscr{P}_{m}(q,\sigma) & =(m-1)\log(1+\frac{t}{2m}(q+\sigma-2)+O(m^{-2}))+\log(1-(t/m)+O(m^{-2}))\\
 & =(m-1)\frac{t}{2m}(q+\sigma-2)+O(m^{-1})
\end{align*}
uniformly in $q$ (and pointwise in $t$). Thus
\[
\mathscr{P}_{m}(q,\sigma)\to\exp\left(\frac{t}{2}(q+\sigma-2)\right).
\]
It follows that
\[
\frac{q+1}{2}\mathscr{P}_{m}(q,1)+\frac{q-1}{2}\mathscr{P}_{m}(q,-1)\to e^{t(q-1)/2}\left(\frac{1-e^{-t}}{2}+\frac{1+e^{-t}}{2}q\right)=F_{t}(q).
\]
Thus given $\epsilon>0$ we can choose $m_{0}$ so that for any $m>m_{0}$,
for all $q\in[0,1]$, 
\[
Q_{q}(\theta_{m,t}^{-1}(B))\geq F_{t}(q)-\epsilon,
\]
as claimed.
\end{proof}

\subsection{Verification of the hypothesis of Proposition \ref{prop:TandTheta}.}
\begin{lem}
\label{lem:thetammtworkslocally.}For every $t>0$, every finite tuple
$x_{1},\ldots,x_{n}\in L^{2}(M)$, and every $\varepsilon,\eta>0$,
there exists an $m_{0}$ so that for all $m>m_{0}$
\[
\theta_{m,t}(w)=w,\qquad\|\theta_{m,t}(A)-A\|_{2}<\eta,
\]
and 
\[
\sum_{i=1}^{n}\mathcal{Q}_{q}\bigl(\theta_{m,t}^{-1}(x_{i})\bigr)\ge\sum_{i=1}^{n}\mathcal{Q}_{F_{t}(q)}(x_{i})-\varepsilon\qquad(0\le q\le1).
\]
\end{lem}
\begin{proof}
Put $\theta_{m}=\theta_{m,t}^{-1}$. Fix $g\in\mathbb{F}_{2}$. Then
by Corollary \ref{cor:lowerEstimatePsigVsPsiB}, there is some $m_{0}$
so that for all $m>m_{0}$
\begin{align*}
Q_{q}(\theta_{m}(g)) & \geq Q_{q}(\theta_{m}(B))^{\ell(g)}-C_{g}\eta_{f_{m,t}}\geq Q_{q}(\theta_{m}(B))^{\ell(g)}-C_{g}t/m.
\end{align*}
Thus for any $\epsilon'>0$ we can use Lemma \ref{lem:F_tlowerEqThetaOfB}
to deduce that for $m$ large enough,
\[
Q_{q}(\theta_{m}(g))\geq F_{t}(q)^{\ell(g)}-\epsilon',\qquad\forall q\in[0,1].
\]
If $g\neq h$, applying Lemma \ref{lem:OffDiagonalEstimate} to $R_{m}(g)$,
$R_{m}(h)$ and $\bar{f}_{m,t}$, and noting that $R_{-m}$ is an
isometry and $R_{m}T_{q}=T_{q}R_{m}$ gives
\[
\left|\langle\theta_{m}(g),T_{q}\theta_{m}(h)\rangle\right|\leq2\sqrt{\eta_{\bar{f}}}+\eta_{\bar{f}}\leq2\sqrt{t/m}+t/m
\]
and thus for $m$ large enough,
\[
\left|\langle\theta_{m}(g),T_{q}\theta_{m}(h)\rangle\right|\leq\epsilon',\qquad\forall q\in[0,1].
\]

Now, if $x=\sum\alpha_{g}g\in L^{2}(M)$ is a finite sum, choose $\epsilon'=\epsilon\left[\sum_{g}|\alpha_{g}|\right]^{-2}$
and apply the definition of $Q$ to deduce
\begin{align*}
Q_{q}(\theta_{m}(x)) & \geq\left[\sum_{g}|\alpha_{g}|^{2}F_{t}(q)^{\ell(g)}\right]-\epsilon'\sum_{g}|\alpha_{g}|^{2}-\epsilon'\sum_{g\neq h}|\alpha_{g}||\alpha_{h}|\\
 & =Q_{F_{t}(q)}(x)-\epsilon'\left(\sum_{g}|\alpha_{g}|\right)^{2}=Q_{F_{t}(q)}-\epsilon.
\end{align*}
For a general $x\in L^{2}(M)$, choose a finite sum $x'$ approximating
$x$, and note that since $T_{q}\circ\theta_{m}$ and $T_{F_{t}(q)}$
are both contractions, $|Q_{q}(\theta_{m}(x))-Q_{q}(\theta_{m}(x'))|\leq(\Vert x\Vert_{2}+\Vert x'\Vert_{2})\Vert x-x'\Vert_{2}$
similarly for $|Q_{F_{t}(q)}(x)-Q_{F_{t}(q)}(x')|$. Choosing an even
bigger $m_{0}$ ensures that $\|\theta_{m,t}(A)-A\|_{2}<\eta$. 
\end{proof}
\begin{thm}
\label{thm:freeComplementation}Let $\eta>0$ be given. Let $A,B$
be two free Haar unitaries generating $L(\mathbb{F}_{2})$, and let
$w=ABA^{-1}B^{-1}$, and let $\mathscr{A}=W^{*}(w)$. Then there exists
a Haar unitary $v\in L(\mathbb{F}_{2})$ satisfying $\Vert v-A\Vert_{2}<\eta$,
and so that $L(\mathbb{F}_{2})=\mathscr{A}*W^{*}(v).$
\end{thm}
\begin{proof}
Define recursively $q_{0}=0$, $t_{j}=(1-q_{j})^{2}/4$, $q_{j+1}=F_{t_{j}}(q_{j})$.
Let $\mathbf{T}_{k}=T_{q_{k}}$, $\mathbf{Q}_{k}=Q_{q_{k}}$. We claim
that with this choice of $\mathbf{T}_{k}$ , the hypothesis of Proposition
\ref{prop:TandTheta} is then satisfied.

We first show that $q_{j}$ increase to $1$. Indeed, $F_{0}(q)=q$,
$\partial_{t}F_{t}(q)\Big|_{t=0}=(1-q)^{2}/2$ and $|\partial_{t}^{2}F_{t}(q)|\leq2$
for all $t\geq0$ and $0\leq q\leq1$. Thus by Taylor's theorem,
\[
F_{t}(q)\geq q+\frac{t}{2}(1-q)^{2}-t^{2}.
\]
Substituting $t_{j}=(1-q_{j})^{2}/4$ gives
\begin{equation}
q_{j+1}-q_{j}\geq\frac{(1-q_{j})^{4}}{16}.\label{eq:growthofqj}
\end{equation}
Since $F_{t}(q)\leq1$, the sequence $q_{j}$ must have a limit $a\leq1$.
But if $a<1$, (\ref{eq:growthofqj}) gives that $q_{j+1}-q_{j}\geq(1-a)^{4}/16>0$,
which is a contradiction. Thus $q_{j}$ increase to $1$. 

Thus $\mathbf{T}_{k}\to\operatorname{Id}$ strongly; $\mathbf{T}_{0}=E_{N}$
by definition, and each $\mathbf{T}_{k}$ is a positive contraction. 

Fix $k\geq0$. Given $x_{1},\dots,x_{n}\in L^{2}(M)$, $\eta>0$,
$\epsilon>0$ Lemma \ref{lem:thetammtworkslocally.} gives an automorphism
$\theta=\theta_{m,t_{k}}$ that satisfies $\theta(w)=w$, $\Vert\theta(A)-A\Vert_{2}<\eta$,
and 
\begin{align*}
\sum_{i=1}^{n}\mathbf{Q}_{k}(\theta^{-1}(x_{i})) & =\sum_{i=1}^{n}Q_{q_{k}}(\theta^{-1}(x_{i}))\\
 & \geq\sum_{i=1}^{n}Q_{F_{t_{k}}(q_{k})}(x_{i})-\epsilon=\sum_{i=1}^{n}\mathbf{Q}_{k+1}(x_{i})-\epsilon
\end{align*}
as claimed. 
\end{proof}
\begin{cor}
\label{cor:amalgamatedIsFree}With the notation of the previous Theorem,
let $\mathscr{A}=W^{*}(w)\subset L(\mathbb{F}_{2})$. Let $(P,\tau_{P})$
be any tracial von Neumann algebra, and assume that $\mathscr{A}\subset P$
as a unital von Neumann subalgebra so that the restriction of the
group trace on $L(\mathbb{F}_{2})$ to $\mathscr{A}$ equals the restriction
of $\tau_{P}$. Then
\[
(L(\mathbb{F}_{2}),E_{\mathscr{A}}^{L(\mathbb{F}_{2})}*_{\mathscr{A}}(P,E_{\mathscr{A}}^{P})\cong L(\mathbb{Z})*(P,\tau_{P}),
\]
where $E_{\mathscr{A}}^{L(\mathbb{F}_{2})}$ and $E_{\mathscr{A}}^{P}$
denote the trace-preserving conditional expectations. 
\end{cor}
\begin{proof}
Let $v$ be such that $L(\mathbb{F}_{2})=W^{*}(v)*W^{*}(w)$, and
put $\mathscr{B}=W^{*}(v)\cong L(\mathbb{Z})$. Then
\[
L(\mathbb{F}_{2})*_{\mathscr{A}}P=(\mathscr{B}*\mathscr{A})*_{\mathscr{A}}P=\mathscr{B}*P\cong L(\mathbb{Z})*P,
\]
as claimed.
\end{proof}

\section{Further group von Neumann algebra consequences.}

In addition to the consequences for finitely generated discrte non-elementary
subgroups of $PSL_{2}(\mathbb{R})$ that were already pointed out,
Theorem \ref{thm:freeComplementation} and Corollary \ref{cor:amalgamatedIsFree}
have several other applications.

\subsection{Commutator MASA in $L(\mathbb{F}_{2m})$. }
\begin{prop}
Let $a_{1},v_{1},\dots,a_{n},v_{n}$ be genereators of $\mathbb{F}_{n}$
and let $c_{j}=a_{j}v_{j}a_{j}^{-1}b_{j}^{-1}$, $w=c_{1}\cdots c_{n}$.
Then there exists a Haar unitary $v$ so that 
\[
L(\mathbb{F}_{2n})=W^{*}(w)*W^{*}(v,a_{2},b_{2},\dots,a_{n},b_{n})\cong W^{*}(w)*L(\mathbb{F}_{2n-1}).
\]
Thus the product-of-commutators MASA in $L(\mathbb{F}_{2m})$ is also
freely complemented. 
\end{prop}
\begin{proof}
Let $v\in W^{*}(a_{1},b_{1})$ be such that $W^{*}(a_{1},b_{1})=W^{*}(c_{1},v)$.
Then $c_{1},v,a_{2},b_{2},\dots,a_{n},b_{n}$ are freely independent
generators of $L(\mathbb{F}_{2})$. Let $r=c_{2}\cdots c_{n}\in\langle a_{2},b_{2},\dots,a_{n},b_{n}\rangle$.
Then $c_{1}r,v,a_{2},b_{2},\dots,a_{n},b_{n}$ are again freely independent
generators. Since $w=c_{1}r$, the claimed isomorphism follows.
\end{proof}
With the notation of the previous Proposition, let $\mathscr{A}=W^{*}(w)\subset L(\mathbb{F}_{2n})$.
Let $(P,\tau_{P})$ be any tracial von Neumann algebra, and assume
that $\mathscr{A}\subset P$ and $\tau_{L(\mathbb{F}_{2n})}|_{\mathscr{A}}=\tau_{P}|_{\mathscr{A}}$.
Then
\[
(L(\mathbb{F}_{2n}),E_{\mathscr{A}}^{L(\mathbb{F}_{2n})}*_{\mathscr{A}}(P,E_{\mathscr{A}}^{P})\cong L(\mathbb{F}_{2n-1})*(P,\tau_{P}).
\]
In particular, if $\Lambda$ is any group $x\in\Lambda$ an element
of infinite order. Then $L(\Lambda*_{x=w^{-1}}\mathbb{F}_{2m})\cong L(\Lambda)*L(\mathbb{F}_{2m-1}).$
Thus, using \cite{dykema:freeProductHyperfinite}, we obtain:
\begin{prop}
Let $\Lambda$ be any amenable group, and let $x\in\Lambda$ be an
element of infinite order. Then $L(\Lambda*_{x^{-1}=\prod_{j=1}^{n}a_{j}b_{j}a_{j}^{-1}b_{j}^{-1}}\langle a_{1},b_{1},\dots,a_{n},b_{n}\rangle)\cong L(\mathbb{F}_{2n})$.
\end{prop}

\subsection{Non-orientable surfaces of negative curvature.}

If $\mathscr{N}_{h}$ is a closed non-orientable surface of genus
$h\geq3$, then its fundametal group is given by the presentation
\[
\pi_{1}(\mathscr{N}_{h})=\left\langle a,b,x_{1},\ldots,x_{h-2}\ \middle|\ aba^{-1}b^{-1}x_{1}^{2}\cdots x_{h-2}^{2}=1\right\rangle .
\]
Thus $\pi_{1}(\mathscr{N}_{h})=\langle a,b\rangle*_{aba^{-1}b^{-1}=(x_{1}^{2}\cdots x_{h-2}^{2})^{-1}}\langle x_{1},\dots,x_{h-2}\rangle$.
Since $x_{1}^{2}\cdots x_{h-2}^{2}$ has infinite order in $\mathbb{F}_{h-2}$
we can use Corollary \ref{cor:amalgamatedIsFree} to deduce that $L(\mathscr{N}_{h})\cong L(\mathbb{F}_{h-1})$.
Combined with Theorem \ref{thm:subgroupsPSL2} this gives:
\begin{thm}
Let $S$ be a closed connected surface of negative Euler characteristic.
Then $L(\pi_{1}(S))\cong L(\mathbb{F}_{1-\chi(S)})$.
\end{thm}


\begin{thebibliography}{10}
\bibitem{Boutonnet-Popa:freeComplementation}R. Boutonnet and S. Popa,
\textit{Maximal amenable MASAs of radial type in the free group factors},
arXiv:2302.13355, to appear in Proc. AMS.

\bibitem{Bundgaard-Nielsen}S. Bundgaard, J. Nielsen, \textit{On normal
subgroups with finite index in $F$-groups}. Mat. Tidsskr. B. (1951)
56--58

\bibitem{Chau}T.C. Chau, \textit{A note concerning Fox’s paper on
Fenchel’s conjecture}. Proc. Amer. Math. Soc. 88(4), (1983) 584--586.

\bibitem{conleyGaboriauMarksTuckerDrob:treeabilitySurface} C. Conley,
D. Gaboriau, A. Marks, R. Tucker-Drob, \textit{One-ended spanning
subforests and treeability of groups}, Preprint https://arxiv.org/abs/2104.07431v3
(2021, revised 2026).

\bibitem{dykema:interpolated}K. Dykema, \textit{Interpolated Free
Group Factors}, Pacific J. Math. 163, (1994), 123-135.

\bibitem{dykema:freeProductHyperfinite}K. Dykema, \textit{Free Products
of Hyperfinite von Neumann Algebras and Free Dimension}, Duke Math.
J. 69 (1993), 97-119.

\bibitem{fox}R.H. Fox, \textit{On Fenchel’s conjecture about $F$-groups}.
Mat. Tidsskr. B (1952), 61--65 

\bibitem{gaboriau:stronglyTreeable}D. Gaboriau, \textit{Examples
of groups that are measure equivalent to the free group}, Ergodic
Theory Dynam. Systems 25 (2005), no. 6, 1809--1827.

\bibitem{delaHarpe-Voiculescu:fuchsian}P. de la Harpe, D. Voiculescu,
\textit{A problem on the II$_{1}$ Factors of Fuchsian Groups}, Astérisque,
232 (1995), 155-158.

\bibitem{jones:index} V.F.R. Jones,\textit{ Index for Subfactors},
Invent. Math. 72, (1983) 1-25 

\bibitem{Mennicke}J. Mennicke, \textit{Eine Bemerkung über Fuchssche
Gruppen}. Invent Math 2, (1967) 301--305.

\bibitem{Mennicke:correction}J. Mennicke, \textit{Corrigendum zu
``Eine Bemerkung über Fuchssche Gruppen, Inventiones math.2, 301--305
(1967).''} Invent. Math. 6 (1968/69) p. 106.

\bibitem{nielsen}J. Nielsen, \textit{The commutator group of the
free product of cyclic groups}. Mat. Tidsskr. B. (1948), 49--56

\bibitem{pimsner-popa:index}M. Pimsner, S. Popa, \textit{Entropy
and Index for Subfactors,} Ann. Scient. Ec. Norm. Sup., 46 (1986),
57--106.

\bibitem{popa:maxInjective}S. Popa, \textit{Maximal injective subalgebras
in factors associated with free groups}, Advances in Math., 50 (1983),
27-48.

\bibitem{popa:compressible} S. Popa, \textit{Compressible Subalgebras
in II$_{1}$ Factors}, Comm.. Math. Phys. (2026) 407:221 

\bibitem{radulescu:interpolated}F. Radulescu, \textit{Stable equivalence
of the weak closures of free groups convolution algebra}s, Comm. Math.
Phys. 156, (1993), 17-36 

\end{thebibliography}
\end{document}